\input{style/arxiv-ba.cfg}
\documentclass[ba,linksfromyear,preprint]{imsart}
\usepackage{vtexurl}

\pubyear{2015}
\volume{10}
\issue{1}
\firstpage{223}
\lastpage{226}
\doi{10.1214/14-BA935}

\startlocaldefs
\newcommand{\enquote}[1]{``#1''}
\expandafter\ifx\csname natexlab\endcsname\relax\fi

\let\bid@start\@empty
\let\bid@end\@empty

\def\MR@url{http://www.ams.org/mathscinet-getitem?mr=}
\def\MR#1{\href{\MR@url#1}{MR#1}}
\def\BDOI#1{%
\edef\doi@base@i{\doi@base}\def\doi@base{}%
doi:~\doiurl{\doi@base@i#1}}

\appto\bid@start{\def\doi@size{\ttfamily}}
\appto\bid@end{\unskip.}

\define@key{bid}{mr}{\def\bid@mr{\MR{#1}}}
\define@key{bid}{doi}{\def\bid@doi{\BDOI{#1}}}
\define@key{bid}{pubmed}{}
\define@key{bid}{pii}{}
\define@key{bid}{issn}{}

\def\bid#1{%
       \bgroup
       \bid@start
       \let\bid@output\@empty
       \setkeys{bid}{#1}\ignorespaces%
       \ifdefvoid\bid@mr{}{\appto\bid@output{\bid@mr}}%
       \ifdefvoid\bid@doi{}{
         \ifdefempty\bid@output{}{\appto\bid@output{. }}
         \appto\bid@output{\bid@doi}%
       }%
       \bid@output
       \bid@end
       \egroup
}

\endlocaldefs

\begin{document}

\begin{frontmatter}
\title{Comment on Article by Berger, Bernardo, and~Sun\thanksref{T1}}
\runtitle{Comment on Article by Berger, Bernardo, and Sun}

\relateddois{T1}{Main article DOI: \relateddoi[ms=BA915,title={James O. Berger, Jose M. Bernardo, and Dongchu Sun. Overall Objective Priors}]{Related item:}{10.1214/14-BA915}.}

\begin{aug}
\author[a]{\fnms{Siva} \snm{Sivaganesan}\corref{}\ead[label=e1]{sivagas@ucmail.uc.edu}}

\runauthor{S. Sivaganesan}

\address[a]{University of Cincinnati, \printead{e1}}
\end{aug}

\end{frontmatter}


Congratulations to the authors on this important paper that leads the
way in selecting an objective overall prior for estimation. The paper
is very enjoyable to read.

The authors provide three possible approaches one could use to find an
overall objective prior suitable for use when there is interest in
simultaneous estimation of several parameters. They illustrate the
approaches in several examples, and give a comprehensive evaluation of
the resulting priors. The proposed new approaches are very carefully
thought out, and hold much promise for the development of a single
overall objective prior in many more models. This is a very interesting
paper and is likely to, and hopefully will, spur increased research in
this new development to find overall objective priors for estimation.

Selection of good objective priors is very important in the practice of
Bayesian analysis since, often, there is little or no prior information
available for at least some of the parameters, especially in complex
models with large number of parameters. Use of diffuse priors is not
always good or optimal.
The reference prior approach has been very successful in providing a
way to get objective priors for estimation in numerous standard and
non-standard models. It was introduced in \cite{b79} to derive a
non-informative prior for estimation of a scalar parameter. In simple
terms, the reference prior is the prior that maximizes, in an
asymptotic sense, the missing information in a prior measured by the
Kullback--Leibler distance between the prior and the posterior
distribution. The approach gave good priors in the one-parameter case,
but did not easily extend to multi-parameter cases. A series of
influential articles beginning with \citet{bb89,bb92}, and
later by Berger\xch{,}{and} Bernardo and Sun extended the reference prior
approach to multi-parameter problems. and formalized the approach,
e.g., see \citet{bbs09,bbs12}, \citet{sb98}, and \citet{bs08}.
It is reasonable to say that the reference prior approach is the best
formal approach to obtain an objective prior for estimation.

The literature is now filled with reference priors for several standard
and non-standard models, ready for use when objective Bayesian
estimation is desired.
The reference prior approach has often been found to have the virtue
of giving good priors when the conventional choices fail, for example,
due to the behavior of the likelihood in the tail. One case in point is
in spatial modeling, see \cite{bos01}.
How this is achieved seems to be a mystery to me. In this paper too,
for the multlinomial example using the hierarchical approach, the
reference prior for the hyper parameter turned out to be a proper prior
to compensate for the slow decay of likelihood in the tail.
However, one runs into difficulty in implementing the reference prior
approach when there are more than one parameter of interest.
Given a model, there are many reference priors; one prior for each
parameter of interest, or even a set of priors for each parameter of
interest based on different ordering of the rest of the parameters.
These priors can be different for different parameters, requiring a
user to switch priors depending which parameter(s) one is interested in
estimating. This is not convenient to explain or appealing to use in
practice, when one is interested in estimating more than one parameter
and the corresponding reference priors are different for these
parameters. Having a single objective prior for a given model, that
works well for most natural parameters of interest is desirable.
In this paper, the authors have taken up this important task and have
given three possible approaches to get a single common ``Overall
Objective Prior'' for simultaneously estimating several parameters of interest.

First, the authors set out to identify models for which there is a
unique common reference prior for each of the natural parameters in the
model under different orderings of the rest of the parameters. The
authors give a condition on the the Fisher information matrix for such
a single reference prior to exist, and provide examples which show that
such a common reference prior can exist for the natural parameters of
many different models.

The other two approaches provided in the paper constitute interesting
novel ideas and developments, and include the Reference Distance
approach and the Hierarchical Prior approach.

Hierarchical prior approach assumes a priori that the parameters of
interest, $\theta_i$'s, conditionally on a hyper-parameter $a$, have a
joint proper prior, leaving a prior for $a$ to be determined. When this
conditional prior is in a convenient form in relation to the likelihood
such as a conjugate prior so that the marginal likelihood for $a$ can
be computed in closed form, one can obtain the reference prior for $a$,
which is the Jeffreys prior based on the marginal likelihood.
Then the overall objective prior for $\theta_i$'s is the marginal prior
obtained by integrating the conditional prior for $\theta_i$'s with
respect to the reference prior for $a$.

The reference distance approach is relatively more involved. Suppose
that for each of the parameters of interest $\theta_i$, $i=1,\dots,n$,
one can choose a reference prior. Then the reference distance approach
first postulates a joint parametric family of priors for $(\theta
_1,\dots,\theta_n)$, not necessarily proper priors, indexed by a
hyper-parameter $a$. Then the overall prior is that prior in the family
whose marginal posterior distributions of $\theta_i$'s is closest on
average, in terms of expected Kullback--Leibler distance, to the
marginal reference posteriors of $\theta_i$'s.

The two approaches hold much promise in achieving the goal of finding
overall objective priors for various models and parameters of interest.
The hierarchical approach is particularly appealing, because the
resulting prior itself is a reference prior, and it may also be
relatively easy to derive, which can be a big advantage. However, the
assumption of a convenient hierarchical or exchangeable structure for
the joint prior of the parameters of interest is not always tenable. In
comparison, the derivation of the reference distance approach requires
computation of reference priors for each parameter of interest and a
not-always-easy computation to find the optimal value of $a$, and the
resulting overall objective prior is not necessarily a reference prior.

But, the reference distance approach holds an advantage -- it seems in
most cases one can write down a joint prior for the parameters of
interest, indexed by a suitable hyper-parameter $a$ by inspecting the
reference priors associated with each parameter. As always with the
reference prior, once the hard work is done, it is readily available
for use by everyone. It is a pleasant surprise that the reference
distance approach for the normal model (Section 3.2.4) gives a prior
that is the reference prior for the natural parameters. However, in
general the reference distance approach may yield a prior that is
different from any of the reference priors used in the derivation. Such
an an overall objective prior may also turn out not to have good
posterior behavior for some of the parameters of interest. In some
instances, there may be more than one choice for the parametric class,
each leading to different overall objective prior, and one has to make
a determination which one to use. Would the authors comment on this and
whether they have encountered such scenarios?

In light of these comments, the recommendation by the authors to use
the common reference prior or the hierarchical approach first, and if
not successful, to try the reference distance approach is noteworthy.

It is surprising that the reference prior for $a$ in the hierarchical
approach to the multinomial example turns out to be a proper prior,
making up for the behavior of the marginal likelihood being bounded away from 0 at
infinity. As indicated before, the phenomenon that the reference prior
distributes its mass selectively compensating for the the likelihood's
slow decays in some tail regions is indeed amazing. Perhaps, the
authors can give some general insight into this phenomenon. Both
approaches have been illustrated for the multinomial example, yielding
different overall objective priors. The reference distance approach
sets $a=1/m$, and the reference prior for $a$ in the hierarchical
approach also seems to favor small values for $a$ for large $m$.
However, for moderate values of $m$, the uncertainty in $a$ induced by
the hierarchical prior approach would have an influence on the
estimation of the parameters of interest, may be of an adaptive nature,
unlike in the reference distance approach. Can the authors comment on
this and how one may choose between the two choices?

While the hierarchical prior approach has its advantages, it appears
that there may be more than one choice for the joint distribution for
the parameters of interest,
$\theta_i$'s, in terms of the second stage parameter $a$. In such
cases, one would have to determine what would be the best choice. For
example, in the multi-normal example in Section 4.3, one may
alternatively use $\mu_i \stackrel{iid}{\sim} N(\mu_0,\tau^2)$ with
known $\mu_0$, or $N(\mu,\tau_0^2)$ with known $\tau_0$, or more
generally $N(\mu,\tau^2)$. In the case of $N(\mu_0,\tau^2)$, it appears
that the resulting estimates for individual $\mu_i$'s would shrink
towards $\mu_0$. Is there any particular justification for the choice
of $\mu_0=0$?

\nocite{*}

\bibliographystyle{ba}

\end{document}